\documentclass [12pt]{article}
\usepackage {color}

\begin{document}
\begin{center}
\textbf{Analytic continuation of functions along parallel algebraic curves.}
\end{center}

\begin{center}
\textit{S.A.Imomkulov}\footnote{ Navoi State Pedagogical Institute (Uzbekistan). E-mail:
\underline {\textcolor{blue}{sevdi@rambler.ru}}}  \textit{ and J.U.Khujamov}\footnote{
Urgench State University (Uzbekistan). E-mail: \underline
{\textcolor{blue}{jumanazar-1971@mail.ru}}}
\end{center}

\bigskip

\begin{center}
\textbf{Introduction}
\end{center}

Classical lemma of Hartogs (see [1]), about continuation along fixed
direction, conform, that if holomorphic in a domain

\[
D\times V = D\times {\left\{ {w \in C:\;{\left| {w} \right|} < r} \right\}}
\subset C_{z}^{n} \times C_{w}
\]

\noindent
function $f(z,\;w)$, in each fixed $z \in D$ holomorphic in a ball ${\left\{
{w \in C,\;{\left| {w} \right|} < r} \right\}}$, then it holomorphic on a
collection of variables in $D\times V$. The Hartogs' lemma has a plenty
generalizations which different on characterizations and direct verges to
themes, connected with holomorphic continuation on fixed direction. A
further results in this field are given in a works of Rothstein [2], M.V.
Kazaryan [3], A.S.Sadullaev and E.M.Chirka [4], T.T.Tuychiev [5],
S.A.Imomkulov and J.U.Khujamov [6], S.A.Imomkulov [7] and etc.

In this work we study a problem about analytic continuation along parallel
algebraic curves.

Algebraic curve in ${\rm {\bf C}}^{{\rm {\bf 2}}}$ will be determining as a
set of zero of some of polynomials:

\[
A = \;{\left\{ {\left( {\xi ,\eta}  \right) \in {\rm {\bf C}}^{{\rm {\bf
2}}}:\;\;P\left( {\xi ,\eta}  \right) = 0} \right\}}.
\]

A set of regular points $A^{0}$of algebraic curve $A$, open on $A$, and the
set of critical points $A^{c} = A\backslash A^{0}$ is discrete set (see [1],
[8],[9] ).

The algebraic set call irreducible, if it is impossible to present in the
view of joined algebraic sets, which differ from it.

\textbf{Theorem.} Let $D$ is a domain from ${\rm {\bf C}}^{{\rm {\bf n}}}$
and $V$ is a domain from some of irreducible algebraic curve $A$. If a
function $f(z,w)$ holomorphic in a domain $D\times V \subset C_{z}^{n}
\times A$ and in each fixed $\xi $ from some of nonpluripolar set $E \subset
D$, a function $f(\xi ,w)$ variable of $w$ continuous till the function,
holomorhpic on a whole $A$, excepting finite sets of singularities (from
$A^{0})$, then $f(z,w)$ holomorphic extended in $(D\times A)\backslash S$,
where $S$- some of analytic subset $D\times A$.

\begin{center}
\textbf{\S 1. Holomorphic functions on algebraic set.}
\end{center}

\bigskip

Any algebraic curve $A \subset {\rm {\bf C}}^{{\rm {\bf 2}}}$ in
sufficiently small neighborhood $U$ (each its point $a \in A)$ is ramified
covering over neighborhood of projection point $a$ in some plane $C \subset
{\rm {\bf C}}^{{\rm {\bf 2}}}$ (see. [1], chapter 2. \S 8. Page
175).

Observe, that on algebraic curves defining correctly not only holomorphic
functions, its derivatives also. Holomorphic function and its derivatives in
a point  defining with the help following equality:

\[
D^{k}f{\left| {_{a}}  \right.} = {\frac{{\partial ^{k}}}{{\partial
z^{k}}}}f(\pi ^{ - 1}(z)){\left| {_{U}}  \right.},
\quad
k = 0,\,1,\,2,\,3,\,...,
\]

\noindent
where $\pi :\,A \to {\rm {\bf C}}$ locally biholomorphic mapping, called
projection and $U$ -- neighborhood of point $a$, in which a restriction $\pi
^{}{\left| {_{U}}  \right.}$ - biholoorphic, $\pi ^{ - 1}{\left| {_{U}}
\right.}$ - inverse mapping and in second member is a derivative in a point
$z = \pi (a)$.

\begin{center}
\textbf{\S 2. The Jacoby -- Hartogs series}
\end{center}

\bigskip

We consider the function $f(z,w)$ holomorphic in a domain $D\times V,\,\,\,D
\subset C^{n},\,\,V \subset A$. Assume, that $\pi ^{ - 1}(0) \in V$. Let
$g(\zeta )$ - rational function on $\zeta \in C$ such that $g(0) = 0$. Then
at small $\rho $ there exist connected component $\Lambda _{\rho}  $ of set
${\left\{ {\zeta :\,\,{\left| {g(\zeta ) < \rho}  \right|}} \right\}}$ such
that $O \in \Lambda _{\rho}  $, $\pi ^{ - 1}(\Lambda _{\rho}  ) \subset V$.
Since $f(z,w)$ holomorphic in the domain $D\times \pi ^{ - 1}(\Lambda _{\rho
} )$ then, at every fixed point $z \in D$ it can be expanded in a series
Jacoby - Hartogs (see [4])

\begin{equation}
\label{eq1}
f(z,w) = {\sum\limits_{k = 0}^{\infty}  {c_{k} (z,\pi (w))g^{k}(\pi (w))}
}\,,
\end{equation}

\noindent
where an coefficients of series defining as follows

\[
c_{k} (z,w) = {\frac{{1}}{{2\pi i}}}{\int\limits_{\partial (\pi ^{ -
1}(\Lambda _{\rho}  )} {f(z,\,\pi (\eta )){\frac{{g(\pi (\eta )) - g(\pi
(w))}}{{g^{k + 1}(\pi (\eta ))(\eta - w)}}}d\eta}}
\]

, $k = 0,1,2,.....$. \newline

It follows, that $c_{k} (z,w)$ holomorphic functions in $D\times A$. At a
fixed point $z \in D$ the series (\ref{eq1}) converge in a lemniscates ${\left\{
{{\left| {g(\pi (w))} \right|} < R^{(g)}(z)} \right\}}$, where $R^{(g)}(z)$
defining as follows

\[
R^{(g)}(z) = {\frac{{1}}{{\overline {{\mathop {\lim} \limits_{k \to \infty}
}} \sqrt[{k}]{{{\left\| {c_{k} (z,w)} \right\|}_{\pi ^{ - 1}(K)}}} }}}{\rm
.}
\]

Here $K \subset C$- arbitrary nonpolar compact, which does not hold poles
$g$, and the limit in right-hand member of equality does not depend on
choice of such compact. Notice, that the value $R^{(g)}(z)$ is a maximal
radius, for which the function $f$ is holomorphic inside of
lemniscates${\left\{ {{\left| {g(w)} \right|} < R^{(g)}(z)} \right\}}$.

\textbf{Lemma 1.} The Jacoby -- Hartogs series (\ref{eq1}) converge uniformly inside
of open set

\[
G_{g} = {\left\{ {(z,w) \in D\times A\,:\,\,\,{\left| {g(\pi (w))} \right|}
< R_{\ast} ^{(g)} (z)\,\,} \right\}},\,\,z \in D,
\]

\noindent
where $R_{\ast} ^{(g)} (z) = {\mathop {\underline {\lim}}  \limits_{\xi \to
z}} R^{(g)}(\xi )$ - normalization from below. The function $ - \ln R_{\ast
}^{(g)} (z)$ is plurisubharmonic in $D$, $R_{\ast} ^{(g)} (z) \le
R^{(g)}(z),\,\,z \in D$ and the set ${\left\{ {z \in D:\,\,R_{\ast} ^{(g)}
(z) < R^{(g)}(z)} \right\}}$ is pluripolar.

We denote by $\Re = {\left\{ {g(\zeta )} \right\}}$ - countable family of
all rational functions with coefficients from the set $Q + iQ$ ($Q$ - the
set of rational number) such that, each function $\Re = {\left\{ {g(\zeta )}
\right\}}$ has a zero only in a point $w = 0$. In order to study convergence
domain, corresponding Jacoby -- Hartogs series will be useful following
lemma about approximation of flat sets by rational lemniscates.

\textbf{Lemma 2.} ([4],[6],[7]). Let $\Sigma $- close polar set from
$C\backslash {\left\{ {0} \right\}}$ and $K$ - compact in $C\backslash
{\left\{ {\Sigma}  \right\}}$. Then, there exist rational function $g \in
\Re $ such that the lemniscates ${\left\{ {\zeta :\,\,{\left| {g(z)}
\right|} < 1} \right\}}$ is connected, belong to $C\backslash {\left\{
{\Sigma}  \right\}}$ and holding $K$.

\begin{center}
\S \textbf{3. Some properties of pseudoconcave sets.}
\end{center}

\bigskip

The properties of pseudoconcave sets has been studied in works [10-13]. Let
$S$- pseudoconcave subset of domain $D\times V$. Assume, that $S$ does not
cross $D\times \partial V$ and

\[
\,\,S_{a} = S \cap {\left\{ {z = a} \right\}}.
\]

Then:

1) the function $\ln (capS_{z} )$, where cap -- capacity (
\^{\i}\'{a}\^{\i}\c{c}\'{\i}\`{a}$\div $\`{a}{\aa}\`{o}
{\aa}\`{\i}\^{e}\^{\i}\~{n}\`{o}\"{u} (transfinite diameter) of flat set, is
plurisubharmonic in $D$(see [13]).

2) if $S_{z} $ is finite for all $z$ from some pluripolar set $E \subset D$,
then $S$ is analytic set (see[12]).

3) if $S_{z} $ is polar for all $z$ from some nonpluripolar set $E \subset
D$, then $S$ is pluripolar set (see [10], [11]).

\textbf{Defenition.} The close set $S \subset D\times A$ called
pseudoconcave set in the domain $D\times A$ if for any point $a \in S$ there
exist some neighborhood $U \subset D\times A$ and holomorphic in
$U\backslash S$ function $f$ such that it does not converge holomorphically
to the point $a$.

\textbf{Lemma 3}. Let $D$ - a domain from ${\rm {\bf C}}^{{\rm {\bf n}}}$
and $V$ - a domain from irreducible algebraic curve $A$, such that $0 \in
\pi (V)$. Let $S$ - pseudoconcave subset of the domain $D\times V$. Assume,
that $S$ does not cross $D\times \partial V$. Then , if $S_{z} $ - finite
for all $z \in D$, then $S$- analytic subset of the domain $D\times V$.

\begin{center}
\textbf{\S 4. The proof of theorem.}
\end{center}

\bigskip

We expand of function $f(z,w)$ in Jacoby -- Hartogs series by degrees of
function $g(\pi (w))$, ($g(\zeta ) \in \Re \,\,,\,\,\,\zeta = \pi (w),\,\,0
\in \pi (V))$,

\begin{equation}
\label{eq2}
f(z,w) = {\sum\limits_{k = 0}^{\infty}  {c_{k} (z,w)g^{k}(\pi (w))}} \,{\rm
,}
\end{equation}

\noindent
where $c_{k} (z,w) \in O(D\times V)$. It is possible, because $f(z,w)$ -
holomorphic in $D\times V$ and at sufficiently small $\rho > 0$ the
lemniscates ${\left\{ {w:\,\,\,\,\,{\left| {g(w)} \right|} < \rho}
\right\}}$ belong to $V$. According to lemma 1 the series (\ref{eq2}) uniformly
converge inside of set

\[
G_{g} = {\left\{ {(z,w)\,:\,\,\,{\left| {g(\pi (w))} \right|} < R_{\ast
}^{(g)} (z)\,\,} \right\}},\,\,z \in D
\]

\noindent
and consequently, its sum holomorphic in it. According to definition of
family of rational functions $\Re $ the set $G_{g} $ is a domain, which
contain $D\times {\left\{ {\pi ^{ - 1}(0)} \right\}}$. The sum of
constructed series (\ref{eq2}) is coincident with $f(z,w)$ in neighborhood $D\times
{\left\{ {\pi ^{ - 1}(0)} \right\}}$ and, thus, (\ref{eq2}) is holomorphic
continuation of $f(z,w)$ in $G_{g} $.

2. Let $g_{1} \,\,\,\,\,\,$and $g_{2} $an arbitrary rational functions from
family $\Re $ and let $f_{1} (z,w)\,$ and $\,f_{2} (z,w)$ are analytic
continuation of function $f(z,w)$ in a domains $G_{g{}_{1}} $ and $G_{g_{2}
} $correspondingly. Since, for any point $z^{0} \in D$ the function $f_{1}
(z^{0},w)$ single-valued by $w$ and $f_{1} (z^{0},w) = f_{2} (z^{0},w) =
f(z^{0},w)$ for any $(z^{0},w) \in G_{g_{1}}  \cap G_{g_{2}}  $, then
$f(z,w)$ holomorphic in $(G_{g_{1}}  \cup G_{g_{2}}  ) \cap {\left\{ {z =
z^{0}} \right\}}$ so, $f(z,w)$ uniquely converge in $G_{g{}_{1}} \cup
\,\,G_{g_{2}}  $, and it follows that the function $f(z,w)$ uniquely
converge in domain $G = \cup G_{g} $, where union is taken by all rational
functions family $\Re $.

3. Since, at each fixed point $z^{0} \in D$ the function $f(z^{0},w)$
single-valued in $A$, then it follows that analytic continuation of function
$f(z,w)$ (in $D\times A)$ is uniquely. Let $\tilde {G} \subset D\times A$ an
original domain for existence of function $f(z,w)$ regarding to $D\times A$.
So, $\tilde {G}\,\,$ nonexpendable holomorphic in every point $(z^{0},w^{0})
\in S = (D\times A)\backslash \tilde {G}$. From here we receive, that $S$ -
pseudoconcave subset of the domain $D\times A$.

4. Now, using lemma 2, we show, that for any point $z^{0} \in D$, a
set of singular point of function $f(z^{0},w)$ of variable $w \in A$
coincide with layer $S_{z^{0}} $ of the set $S$. In fact, by the
terms of theorem singular set $\Lambda $ of function $f(z^{0},w)$
consist finite number of points, then according to lemma 2, for any
compact $K \subset A\backslash \Lambda $, there exist rational
function $g \in \Re $ such that, the lemniscates ${\left\{ {{\left|
{g(\pi (w))} \right|} < 1} \right\}}$ contains $K$. Consequently,
the lemniscates ${\left\{ {{\left| {g(\pi (w))} \right|} <
R^{(g)}(z^{0})} \right\}}$, so and ${\left\{ {{\left| {g(\pi (w))}
\right|} < R_{\ast} ^{(g)} (z^{0})} \right\}}$ contains (since,
$R_{\ast} ^{(g)} (z^{0})\ge R_{}^{(g)} (z^{0}))$.

5. Let $\Omega $ - an image of domain $\tilde {G} = (D\times A)\backslash S$
in maping $(z,\pi ^{ - 1}(\zeta )) \to (z,\pi ^{ - 1}({\frac{{1}}{{\zeta
}}}))$. The set $(D\times {A}')\backslash \Omega = L$ is also pseudoconcave.
Since, $S$ does not cross the set $D\times {\left\{ {\pi ^{ - 1}(0)}
\right\}}$, then $L$ is bounded and intersection $L \cap {A}'$ for any point
$z^{0} \in D$, consist from finite number of points, i.e. the set $L$ is
satisfying all conditions of lemma 3. Consequently, $L$ - analytic set from
here it is easy to see, that $S$- analytic. The proof of theorem is
complete.

\begin{center}
\textbf{REFERENCES}
\end{center}

\bigskip

Shabat B.V. Introduction to complex analysis. Part 2. Moscow, ``Nauka'',
1985.

Rothstein W. Ein neuer Beweis des hartogsshen hauptsatzes und sline
ausdehnung auf meromorphe functionen // Math. Z. -- 1950.-- V. 53. -- P. 84

Kazaryan M.V. On holomorpgic extension of functions with special
singularities in $C^{n}$. Doc. Acad. Nauk Arm.SSR. 1983. v. 76. p. 13-17.

Sadullaev A.S. and Chirka E.M. On extension of functions with polar
singularities. Math. Sb. 1987. v. 132(174) $^{1}$3. p. 383-390.

Tuychiev T.T. and Imomkulov S.A. Holomorphic extension of functions, having
singularities on parallel multidimensional sections. Doc. Acad. Nauk of
Uzbekistan. 2004. $^{1}$2. p.12-15.

Imomkulov S.A ., Khujamov J. U. On holomorphic continuation of functions
along boundary sections. Mathematica Bohemica. (Czech Republic) -- 2005. V.
130(3) P. 309-322.

Imomkulov S.A. On holomorphic continuation of functions, given on boundary
beam of complex line// Izvestiya Russian Academy of Science. Series of math
-- 2005. -- V. 69, $^{1}$2. -- p.125 -144.

Stoilov S. The theory of functions of complex variables. Volume 1.
Moscow-1962.

Chirka E.M. Complex analytic sets. Moscow-1985.

\end{document}